\documentclass[12pt]{amsart}
\usepackage[latin9]{inputenc}
\usepackage{mathrsfs}
\usepackage{amstext}
\usepackage{amsthm}
\usepackage{amssymb}
\usepackage{wasysym}
\usepackage{esint}

\makeatletter
\numberwithin{equation}{section}
\numberwithin{figure}{section}
\theoremstyle{plain}
\newtheorem{thm}{\protect\theoremname}
  \theoremstyle{plain}
  
  \theoremstyle{plain}

\@ifundefined{date}{}{\date{}}
\usepackage[colorlinks=true, bookmarks=true,pdfstartview=FitV, 
linkcolor=blue, citecolor=blue, urlcolor=blue]{hyperref}

\usepackage{upref}\usepackage{enumerate}\usepackage{tikz}
\usepackage{amscd}\usepackage{amsxtra}\usepackage{latexsym}
\usepackage{yfonts}
\usepackage{mathrsfs}\usepackage{bbm}
\usepackage{bm}
\usepackage{eso-pic}
\usepackage{graphicx}

\def\Xint#1{\mathchoice
{\XXint\displaystyle\textstyle{#1}}%
{\XXint\textstyle\scriptstyle{#1}}%
{\XXint\scriptstyle\scriptscriptstyle{#1}}%
{\XXint\scriptscriptstyle\scriptscriptstyle{#1}}%
\!\int}
\def\XXint#1#2#3{{\setbox0=\hbox{$#1{#2#3}{\int}$}
\vcenter{\hbox{$#2#3$}}\kern-.5\wd0}}

\def\avgint{\Xint-}

\calclayout

\usepackage{pdfsync}

\numberwithin{equation}{section}
\allowdisplaybreaks

\newcommand{\ra}{\rightarrow}

\newcommand{\bey}{\begin{eqnarray*}}
\newcommand{\eey}{\end{eqnarray*}}
\newcommand{\ba}{\begin{align}}
\newcommand{\ea}{\end{align}}
\newcommand{\bea}{\begin{align*}}
\newcommand{\ena}{\end{align*}}
\newcommand{\be}{\begin{equation}}
\newcommand{\ee}{\end{equation}}
\newcommand{\R}{\mathbb R}

\newcommand{\ep}{\epsilon}

\newcommand{\bc}{\begin{center}}
\newcommand{\ec}{\end{center}}

\providecommand{\theoremname}{Theorem}

\makeatother

  \providecommand{\corollaryname}{Corollary}
  \providecommand{\propositionname}{Proposition}
\providecommand{\theoremname}{Theorem}

\begin{document}

\author{Adam Mair }

\address{Adam Mair \\
 Department of Mathematics \\
 University of Alabama \\
 Tuscaloosa, AL 35487, USA}

\email{acmair@crimson.ua.edu}

\author{Kabe Moen}

\address{Kabe Moen \\
 Department of Mathematics \\
 University of Alabama \\
 Tuscaloosa, AL 35487, USA}

\email{kabe.moen@ua.edu}


\title{Two weight bump conditions for compactness of commutators}

\begin{abstract}
We prove certain two weight bump conditions are sufficient for the compactness of the commutator $[b,T]$ where $b\in CMO$ and $T$ is a Calder\'on-Zygmund operator. This is the first result for compactness in the two weight setting without additional assumptions on the individual weights.
\end{abstract}

\thanks{The second author is supported by Simons Collaboration Grant for
  Mathematicians, 160427. A version of this paper will serve as a chapter in
  the second author's Ph.D. dissertation. }
\maketitle

\section{Introduction}
Let $T$ be a Calder\'on-Zygmund operator (CZO) and $b$ be a $BMO$ function. Define the commutator of $b$ and $T$ to be
$$[b,T]f=bTf-T(bf).$$
In this note we are interested in conditions on $b$ and $(u,v)$ that ensure compactness of $[b,T]$ between two different weighted spaces $L^p(v)$ and $L^p(u)$. 

The boundedness of the commutator is well established. Coifmann, Rochberg, and Weiss \cite{CRW} showed that $[b,T]$ is bounded on $L^p(\R^n)$ when $b\in BMO$.  Chung, Pereyra, and P\'erez \cite{MR2869172} showed that $[b,T]$ is bouneded on $L^p(w)$ when $w$ belongs to the Muckenhoupt $A_p$ class and gave the sharp quantitative control of the constants. The two weight boundedness of $[b,T]$ from $L^p(v)$ to $L^p(u)$ is significantly more difficult. Cruz-Uribe and the second author \cite{MR2918187} showed that if one strengthens the two weight $A_p$ condition with certain logarithmic bumps then $[b,T]$ is a bounded linear operator from $L^p(v)$ to $L^p(u)$ (see below for relevant definitions).  We also refer readers to more recent results by Lerner, Ombrosi, and Riveria-R\'ios \cite{LORR}, who improved the results in \cite{MR2918187}.  Finally, we mention in passing the work of  Isralowitz, Pott, and Treil \cite{IsPoTr} and Cruz-Uribe, the second author, and Tran \cite{CMT} in which the two weight boundedness of $[b,T]$ and higher order commutators is obtained assuming more general assumptions on $b$ and the weights $(u,v)$.  

Recall that a linear operator between two Banach spaces $T:X\ra Y$ is compact if $T(B_X)$ has compact closure in $Y$ (here $B_X$ is the unit ball in $X$). The compactness of commutators was first investigated by Uchiyama \cite{MR467384}, who was motivated by applications to complex analysis and Hankel operators.  Uchiyama showed that if $b\in CMO$, the closure of $C_c^\infty(\R^n)$ in the $BMO$ norm, then $[b,T]$ is a compact operator on $L^p(\R^n)$. The original statement in \cite{MR467384} was that $b\in VMO$, the space \emph{vanishing mean oscillation}, is sufficient for compactness of $[b,T]$ on $L^p(\R^n)$. However the spaces $CMO$ and $VMO$ do not coincide over $\R^n$ and the correct space for compactness of commutators is $CMO$. This misnomer has appeared several times in the literature (including \cite{ClCr}) as the article by Bourdaud \cite{MR1927078} points out. Wu and Yang \cite{MR3834654} show that $CMO$ actually characterizes the compactness of $[b,T]$ on $L^p(\R^n)$ for specific operators $T$ like the Riesz transforms. Iwaniec \cite{MR1187088} proved a Fredholm alternative for Beltrami equations as an application of the compactness of $[b,T]$.  Cruz and Clop \cite{ClCr} were also interested in applications to quasiregular mappings and showed that commutator $[b,T]$ is compact on $L^p(w)$ when $b\in CMO$ and $w\in A_p$.  Recently Hyt\"onen and Lappas \cite{MR4383118} have shown an extrapolation theorem for compact operators on weighted spaces. Finally we mention some very recent two weight results which are related to the so-called Bloom $BMO$.  Namely, under certain weighted $BMO$ assumptions on $b$, it is shown that the operator $[b,T]$ is a compact operator from $L^p(v)$ to $L^p(u)$ provided one assumes the side condition $u,v\in A_p$.  We refer the interested reader to Chen, Lacey, Li, and Vempati \cite{ChenLaceyLiVempati} and Hyt\"onen, Oikari, and Sinko \cite{HytonenOikariSinko} for more on compactness in this setting.  Our main result is the following. 

\begin{thm}\label{main} Suppose $T$ is a Calder\'on-Zygmund operator and $b\in CMO$.  If $u$ and $v$ are weights that satisfy
\begin{equation}\label{bump} \sup_Q \|u^{\frac1p}\|_{L^p(\log L)^{2p-1+\delta},Q} \|v^{-\frac1p}\|_{L^{p'}(\log L)^{2p'-1+\delta},Q}<\infty\end{equation}
for some $\delta>0$, then the operator $[b,T]$ is a compact operator from $L^p(v)$ to $L^p(u)$.
\end{thm}

Condition \eqref{bump} is a strengthening of the classic two weight $A_p$ condition.  It does not require individual assumptions on $u$ and $v$.  For example, $u$ can be \emph{any} weight and $v=M^k u$ where $k=\lfloor 2p \rfloor+1$ and $M^k$ is the $k$-fold composition of the Hardy-Littlewood maximal operator $M$ (see \cite{MR1481632}).  We also remark that conditions for compactness of the higher order iterated commutators are more in-depth and will be addressed in a forthcoming paper.

\section{Preliminaries} We begin with the space of functions of bounded mean oscillation. A measurable function $b$ is said to belong to $BMO$ if
$$\|b\|_{BMO}=\sup_Q\avgint_Q \Big|b(x)-\avgint_Q b\Big|\,dx<\infty.$$
The quantity $\|b\|_{BMO}$ is a semi-norm which can be made into a norm by considering $BMO$ modulo constants. Henceforth, we will just refer to $\|b\|_{BMO}$ as the norm of $b$ in $BMO$. The space $CMO$ is the subspace obtained by taking the closure of $C^\infty_c(\R^n)$ with respect to the norm $\|\cdot\|_{BMO}$.

The Hardy-Littlewood maximal operator given by 
$$Mf(x)=\sup_{Q\ni x}\avgint_Q |f(y)|\,dy.$$
A Calder\'on-Zygmund operator is an $L^2(\R^n)$ bounded operator associated to a kernel in the sense
$$Tf(x)=\int_{\R^n} K(x,y)f(y)\,dy, \quad f\in L^\infty_c(\R^n), x \notin \mathsf{supp}\,f.$$
We will assume the kernel $K(x,y)$, defined on $\{(x,y):x\not=y\}$, satisfies the size estimate 
$$|K(x,y)|\leq \frac{C}{|x-y|^n}$$
and smoothness estimate
$$|\nabla K(x,y)|\leq \frac{C}{|x-y|^{n+1}}.$$
We will also need the maximal truncation operator which is given by
$$T^\sharp f(x)=\sup_{\eta>0}\left|\int_{|x-y|>\eta} K(x,y)f(y)\,dy\right|.$$
It is well-known that the Hardy-Littlewood maximal operator, CZOs, and commutators of CZOs with $BMO$ functions are bounded on $L^p(\R^n)$ for all $1<p<\infty$. Moreover, these operators are also bounded on $L^p(w)$ when the weight $w$ belongs to $A_p$:
$$\sup_Q\left(\avgint_Q w\right)\left(\avgint_Q w^{1-p'}\right)^{p-1}<\infty.$$

Following the book \cite{MR2797562}, we will work with pairs of weights $(u,v)$, that is, non-negative $L^1_{\mathsf{loc}}(\R^n)$ functions, such that $u>0$ on a set of positive measure and $v>0$ almost everywhere. We will work a strengthened version of the two weight $A_p$ condition
$$\sup_Q\left(\avgint_Q u\right)\left(\avgint_Q v^{1-p'}\right)^{p-1}<\infty.$$
In order to understand the bump conditions we need some background on Young functions and Orlicz spaces.  A Young function, $\Phi:[0,\infty)\ra [0,\infty)$, is an increasing, convex function that satisfies $\Phi(0)=0$ and $\Phi(t)/t\ra \infty$ as $t\ra \infty$.  Given a Young function $\Phi$ one may define the Orlicz average of $f$ over a cube $Q$
$$\|f\|_{\Phi,Q}=\inf\left\{\lambda>0:\avgint_Q \Phi\Big(\frac{|f(x)|}{\lambda}\Big)\,dx\leq 1\right\}$$
and this quantity is a norm.  For our purposes we will only need Young functions of the form $\Phi(t)=t^p\big[\log(e+t)\big]^a$, where $p>1$ and $a\in \R$. In this case we will use the notation
$$\|f\|_{\Phi,Q}=\|f\|_{L^p(\log L)^a,Q}$$
and refer to such averages as logarithmic bumps.  

In P\'erez \cite{perez94}, it is shown that if the pair of weights satisfies the condition
\begin{equation}\label{maxbump}\sup_Q\left(\avgint_Q u\right)\|v^{-\frac1p}\|_{L^{p'}(\log L)^{p'-1+\delta},Q}<\infty,\end{equation}
for some $\delta>0$ then $M:L^p(v)\ra L^p(u)$.  For CZOs bump conditions on both weights are needed. Lerner \cite{Lern2012} showed that the condition
\begin{equation}\label{czobump}\sup_Q\|u^{\frac1p}\|_{L^p(\log L)^{p-1+\delta},Q}\|v^{-\frac1p}\|_{L^{p'}(\log L)^{p'-1+\delta},Q}<\infty\end{equation}
for some $\delta>0$ is sufficient for the boundedness of $T$ and $T^\sharp$ from $L^p(v)$ to $L^p(u)$.  Finally, for commutators of $BMO$ functions the condition
\begin{equation}\label{commbump}\sup_Q\|u^{\frac1p}\|_{L^p(\log L)^{2p-1+\delta},Q}\|v^{-\frac1p}\|_{L^{p'}(\log L)^{2p'-1+\delta},Q}<\infty\end{equation}
is sufficient for the boundedness of $[b,T]$ from $L^p(v)$ to $L^p(u)$ (see \cite[Theorem 1.3]{MR2918187}).  We note that the condition \eqref{commbump} is stronger than both \eqref{czobump} and \eqref{maxbump}. In particular, if $(u,v)$ satisfy \eqref{commbump}, then all of the operators $M,T,$ and $T^\sharp,$ are bounded from $L^p(v)$ to $L^p(u)$ as well as the commutator $[b,T]$ when $b\in BMO$. We will use this fact in our proof.  We also point out that this is the benefit of the bump conditions, which are universal conditions that are not operator specific.

To prove compactness of the operator $[b,T]$ we will need a characterization of compact subsets of Lebesgue spaces, known as the Kolomogorov-Riesz theorem.  The survey paper by Hanch-Olsen and Holden \cite{MR2734454} contains the basic statement and history of the original theorem.  A characterization of compact subsets of weighted Lebesgue spaces has a more complicated development and was just recently resolved. Early weighted versions of the Kolomogorv-Riesz theorem required extra assumptions on the weight $u$, such as $u\in A_p$, $u\in A_\infty$, or $u^{1-p_0'}\in L^1_{\textsf{loc}}(\R^n)$ for some $p_0$ (see \cite{ClCr} and \cite{MR4303944}).  Guo and Zhao \cite{MR4108843} proved a general Kolomogorov-Riesz theorem in Banach function spaces that generalize $L^p(u)$.  This generalization is paramount for proving two weight compactness, because we are not assuming additional assumptions on the weight $u$.

\begin{thm}[\cite{MR4108843}] \label{KolRiesz} Let $1\leq p<\infty$, $u$ be a weight, and $\mathcal F\subseteq L^p(u)$.  If the family $\mathcal F$ satisfies 
\begin{enumerate}
\item $\mathcal F$ is a bounded subset, i.e.,
$$\sup_{f\in \mathcal F} \|f\|_{L^p(u)}\leq C;$$
\item $\mathcal F$ uniformly vanishes at infinity, that is,
$$\lim_{N\ra \infty}\sup_{f\in \mathcal F}\int_{|x|>N} |f(x)|^p u(x)\,dx=0;$$
\item $\mathcal F$ is equicontinuous, that is, 
$$\lim_{h\ra 0}\sup_{f\in \mathcal F}\int_{\R^n} |f(x+h)-f(x)|^pu(x)\,dx<\ep,$$
\end{enumerate}
then the family $\mathcal F$ has compact closure in $L^p(u)$.
\end{thm} 

Finally, it will be useful to regularize our operators, an idea that goes back to Krantz and Li \cite{MR1835563} (see also \cite{MR3326579}). Given $\eta>0$, we let $K_\eta(x,y)$ be a smooth truncation of the $K(x,y)$ such that 
\begin{enumerate}
\item $K_\eta(x,y)=0$ if $|x-y|\leq \eta$;
\item $K_\eta(x,y)=K(x,y)$ if $|x-y|> 2\eta$;
\item $K_\eta$ satisfies the same size and regularity estimates as $K$, namely
\begin{equation}\label{sizesmooth} |K_\eta(x,y)|\leq \frac{C}{|x-y|^n} \quad \text{and} \quad |\nabla K_\eta(x,y)|\leq \frac{C}{|x-y|^{n+1}}.\end{equation}
\end{enumerate}
We let $T^\eta$ be the operator associated to the kernel $K_\eta$.  We also note that $T^\eta$ is a bounded operator from $L^p(v)$ to $L^p(u)$ whenever $M$ and $T^\sharp$ are also bounded, for example if $(u,v)$ satisfy \eqref{czobump} or \eqref{commbump}.  Indeed,
\begin{multline*}|T^\eta f(x)|\leq \left|\int_{\eta<|x-y|\leq 2\eta} K_\eta(x,y)f(y)\,dy\right|+\left|\int_{2\eta<|x-y|} K(x,y)f(y)\,dy\right|\\
\leq \int_{\eta<|x-y|\leq 2\eta}\frac{|f(y)|}{|x-y|^n}\,dy+T^\sharp f(x)\leq \frac{1}{\eta^n}\int_{|x-y|\leq 2\eta}|f(y)|\,dy+T^\sharp f(x)\\ \leq CMf(x)+T^\sharp f(x).\end{multline*}

\section{Proof of Theorem \ref{main}}
We will make several reductions to prove the compactness of $[b,T]$ from $L^p(v)$ to $L^p(u)$.  First we use the well-known fact that the space of compact operators, $\mathcal K(L^p(v),L^p(u))$, is a closed subset of the space of bounded operators, $\mathcal B(L^p(v),L^p(u))$, in the operator norm topology. The precise bound for $[b,T]$ from \cite{MR2918187} is that if $b\in BMO$ and $(u,v)$ satisfy \ref{commbump}, then
$$\|[b,T]f\|_{L^p(u)}\leq C\|b\|_{BMO}\|f\|_{L^p(v)}$$
and in particular
$$\|[b,T]\|_{\mathcal B(L^p(v),L^p(u))}\leq C\|b\|_{BMO}.$$
The commutator is linear in $b$ which implies
$$\|[b_1,T]-[b_2,T]\|_{\mathcal B(L^p(v),L^p(u))}=\|[b_1-b_2,T]\|_{\mathcal B(L^p(v),L^p(u))}\leq C\|b_1-b_2\|_{BMO}.$$
Thus if $b\in CMO$ and $\{b_j\}$ is a sequence of smooth functions converging to $b$ in $BMO$, then $[b_j,T]\ra [b,T]$ in the operator topology. It is this estimate that breaks down for the higher order iterated commutators.  Next, notice that if $b\in C_c^\infty(\R^n)$ then by using the mean value theorem on $b$ we have
\begin{align*} |[b,T]f(x)-[b,T^\eta]f(x)|&\leq C\int_{|x-y|\leq 2\eta}|b(x)-b(y)|\frac{|f(y)|}{|x-y|^n}\,dy\\
&\leq C \|\nabla b\|_\infty\int_{|x-y|\leq 2\eta}\frac{|f(y)|}{|x-y|^{n-1}}\,dy\\
&\leq C \|\nabla b\|_\infty\sum_{k=0}^\infty \int_{2^{-k}\eta<|x-y|\leq 2^{-k+1}\eta}\frac{|f(y)|}{|x-y|^{n-1}}\,dy\\
&\leq  \eta C\|\nabla b\|_\infty\sum_{k=0}^\infty 2^{-k}\avgint_{B(x,2^{-k+1}\eta)}{|f(y)|}\,dy\\
&\leq \eta C\|\nabla b\|_{\infty}  Mf(x).\\
\end{align*}
This implies that 
$$\|[b,T^\eta]-[b,T]\|_{\mathcal B(L^p(v),L^p(u))}\leq \eta C_b \|M\|_{\mathcal B(L^p(v),L^p(u))}$$
and hence $[b,T^\eta]$ converges to $[b,T]$ as $\eta\ra 0$ in the operator topology. It suffices to show that for a fixed $\eta>0$ and a fixed $b\in C_c^\infty(\R^n)$ that $[b,T^\eta]$ is a compact operator from $L^p(v)$ to $L^p(u)$. Fix $b\in C_c^\infty(\R^n)$ and $\eta>0$ and note that our estimates may depend on $b$ and $\eta$. Consider the unit ball in $L^p(v)$,  
$$B_{L^p(v)}=\{f\in L^p(v):\|f\|_{L^p(v)}\leq 1\}.$$ 
To show that $[b,T^\eta]$ is a compact operator we will show that $\mathcal F=[b,T^\eta](B_{L^p(v)})$ satisfies (a), (b), and (c) of Theorem \ref{KolRiesz}, hence proving that $\mathcal F$ has compact closure in $L^p(u)$.  

Since the operator $[b,T^\eta]$ is bounded from $L^p(v)$ to $L^p(u)$, we have that
$$\|[b,T^\eta]f\|_{L^p(u)}\leq C\|f\|_{L^p(v)}\leq C$$
for $f\in B_{L^p(v)}$ and thus $\mathcal F$ is a bounded subset of $L^p(u)$, i.e., it satisfies (a).  To see that $\mathcal F$ satisfies the uniform vanishing at infinity (b), let $f\in B_{L^p(v)}$ and $Q$ be a cube containing $\mathsf{supp}\,b$. Choose $N>1$ such that $$N>2\sup\{|y|:y\in Q\}.$$  
For such $N$ and $|x|>N$ we have
\begin{align*}||[b,T^\eta]f(x)|&=\left|\int_{\R^n} (b(x)-b(y))K_\eta(x,y)f(y)\,dy\right|\\
&\leq C\|b\|_\infty\int_{\mathsf{supp}\,b}\frac{|f(y)|}{|x-y|^n}\,dy\\ 
&\leq \frac{C\|b\|_\infty}{|x|^n}\int_{Q}|f(y)|\,dy\\
&\leq \frac{C\|b\|_\infty}{|x|^n}\|f\|_{L^p(v)} \left(\int_{Q} v^{-\frac{p'}{p}}\right)^{\frac1{p'}}\\
&\leq \frac{C_{b,v}}{|x|^n}.
\end{align*}
We note that 
$$\left(\int_{Q} v^{-\frac{p'}{p}}\right)^{\frac1{p'}}\leq C|Q|^{\frac1{p'}}\|v^{-\frac1{p}}\|_{L^{p'}(\log L)^{2p'-1+\delta},Q}<\infty$$ 
and hence we have
$$|[b,T^\eta]f(x)|\leq\frac{C_{b,v}}{|x|^{n}}$$
with a finite constant $C_{b,v}$ that depends on $b$ and $v$, but not $f$.  Integrating, over $|x|>N$ yields
$$\int_{|x|>N}|[b,T^\eta]f(x)|^pu(x)\,dx\leq C_{b,v} \int_{|x|>N}\frac{u(x)}{|x|^{np}}\,dx.$$
Since $M(\chi_{[-1,1]^n})(x)\geq C(1+|x|)^{-n}$ and $M$ is bounded from $L^p(v)$ to $L^p(u)$ we have that 
$$\int_{\R^n}\frac{u(x)}{(1+|x|)^{np}}\,dx\leq C\int_{\R^n} M(\chi_{[-1,1]^n})(x)^pu(x)\,dx\leq Cv([-1,1]^n)<\infty.$$
This implies that 
$$\lim_{N\ra \infty}\int_{|x|>N}\frac{u(x)}{|x|^{np}}\,dx=0$$
and thus the family $\mathcal F$ satisfies (b) of Theorem \ref{KolRiesz}.  

We now prove the most difficult of the estimates, the equicontinuity estimate, (c), from Theorem \ref{KolRiesz}.  Let $h\in \R^n$ be sufficiently small to be chosen later and $f\in B_{L^p(v)}$. Consider 
\begin{align*}T^\eta f(x+h)-T^\eta f(x)&=\int_{\R^n}(b(x+h)-b(y))K_\eta(x+h,y)f(y)\,dy\\
&\qquad -\int_{\R^n}(b(x)-b(y))K_\eta(x,y)f(y)\,dy\\
&=(b(x+h)-b(x))\int_{\R^n}K_\eta(x,y)f(y)\,dy \\
&\qquad + \int_{\R^n}(b(y)-b(x+h))(K_\eta(x,y)-K_\eta(x+h,y))f(y)\,dy \\
&=Af(x)+Bf(x),
\end{align*}
where we have added and subtracted the term $b(x+h)T^\eta f(x)$.  The first term is easy to handle with the smoothness of $b$:
$$|Af(x)|\leq |h| \|\nabla b\|_\infty |T^\eta f(x)|\leq C|h|\|\nabla b\|_\infty (Mf(x)+T^\sharp f(x)).$$
Hence 
$$\sup_{\|f\|_{L^p(v)}\leq 1}\|Af\|_{L^p(u)}\leq C|h| $$
which will go to zero uniformly as $h\ra 0$.  For the other term, notice that if $|h|<\frac{\eta}{4}$ and $|x-y|<\frac\eta2$ then $|x+h-y|<|x-y|+|h|<\eta$.  Hence $K_\eta(x+h,y)=K_\eta(x,y)=0$ for such $|h|<\frac\eta4$, and we have 
\begin{align*}
|Bf(x)|&=\left|\int_{|x-y|\geq \frac{\eta}{2}}(b(y)-b(x+h))(K_\eta(x,y)-K_\eta(x+h,y))f(y)\,dy\right|\\
&\leq 2\|b\|_\infty \int_{|x-y|\geq \frac{\eta}{2}}|K_\eta(x+h,y)-K_\eta(x,y)||f(y)|\,dy\\
&\leq 2C|h|\|b\|_\infty \int_{|x-y|\geq \frac{\eta}{2}} \frac{|f(y)|}{|x-y|^{n+1}}\,dy
\end{align*}
where we have used the well-known fact that the gradient estimate \eqref{sizesmooth} implies the standard kernel estimate
$$|K_\eta(x+h,y)-K_\eta(x,y)|\leq C\frac{|h|}{|x-y|^{n+1}}$$
when $|x-y|\geq 2|h|$. Notice $|x-y|\geq 2|h|$ is satisfied since $|x-y|\geq \frac{\eta}{2}\geq 2|h|$.  We now have
\begin{align*} 
\int_{|x-y|\geq \frac{\eta}{2}} \frac{|f(y)|}{|x-y|^{n+1}}\,dy&\leq \sum_{k=0}^\infty \int_{2^{k-1}\eta  \leq |x-y|< 2^k\eta} \frac{|f(y)|}{|x-y|^{n+1}}\,dy\\
&\leq\frac{C}{\eta} \sum_{k=0}^\infty 2^{-k} \avgint_{ B(x,2^k\eta)}|f(y)|\,dy\\
&\leq\frac{C}{\eta}Mf(x).
\end{align*}
The boundedness of the maximal function implies
$$\sup_{\|f\|_{L^p(v)}\leq 1}\|Bf\|_{L^p(u)}\leq C|h|. $$
These estimates do not depend on $f$ and go to zero uniformly as $h\ra 0$ hence we have
$$\lim_{h\ra \infty}\sup_{\|f\|_{L^p(v)}\leq 1}\|[b,T^\eta]f(\cdot+h)-[b,T^\eta]f\|_{L^p(u)}=0.$$
This concludes our proof thus showing that $[b,T]$ is a compact operator.

\bibliographystyle{plain}
\bibliography{CompactnessofCommutators}

\begin{thebibliography}{10}

\bibitem{MR3326579}
\'{A}rp\'{a}d B\'{e}nyi, Wendol\'{\i}n Dami\'{a}n, Kabe Moen, and Rodolfo~H.
  Torres.
\newblock Compact bilinear commutators: the weighted case.
\newblock {\em Michigan Math. J.}, 64(1):39--51, 2015.

\bibitem{MR1927078}
G\'{e}rard Bourdaud.
\newblock Remarques sur certains sous-espaces de {${\rm BMO}(\Bbb R^n)$} et de
  {${\rm bmo}(\Bbb R^n)$}.
\newblock {\em Ann. Inst. Fourier (Grenoble)}, 52(4):1187--1218, 2002.

\bibitem{ChenLaceyLiVempati}
Peng Chen, Michael Lacey, Ji~Li, and Manasa~N. Vempati.
\newblock Compactness of the {B}loom sparse operators and applications.
\newblock 2022.
\newblock arXiv:2204.11990.

\bibitem{MR2869172}
Daewon Chung, M.Cristina Pereyra, and Carlos P\'erez.
\newblock Sharp bounds for general commutators on weighted {L}ebesgue spaces.
\newblock {\em Trans. Amer. Math. Soc.}, 364(3):1163--1177, 2012.

\bibitem{ClCr}
Albert Clop and Victor Cruz.
\newblock Weighted estimates for {B}eltrami equations.
\newblock {\em Ann. Acad. Sci. Fenn. Math.}, 38(1):91--113, 2013.

\bibitem{CRW}
Ronald~R. Coifman, Richard Rochberg, and Guido Weiss.
\newblock Factorization theorems for {H}ardy spaces in several variables.
\newblock {\em Ann. of Math. (2)}, 103(3):611--635, 1976.

\bibitem{MR2797562}
David Cruz-Uribe, Jos\'{e}~Maria Martell, and Carlos P\'{e}rez.
\newblock {\em Weights, extrapolation and the theory of {R}ubio de {F}rancia},
  volume 215 of {\em Operator Theory: Advances and Applications}.
\newblock Birkh\"{a}user/Springer Basel AG, Basel, 2011.

\bibitem{MR2918187}
David Cruz-Uribe and Kabe Moen.
\newblock Sharp norm inequalities for commutators of classical operators.
\newblock {\em Publ. Mat.}, 56(1):147--190, 2012.

\bibitem{CMT}
David Cruz-Uribe, Kabe Moen, and Quan Tran.
\newblock New oscillation classes and two weight bump conditions for
  commutators.
\newblock {\em to appear Collect. Mat.}

\bibitem{MR4108843}
Weichao Guo and Guoping Zhao.
\newblock On relatively compact sets in quasi-{B}anach function spaces.
\newblock {\em Proc. Amer. Math. Soc.}, 148(8):3359--3373, 2020.

\bibitem{MR2734454}
Harald Hanche-Olsen and Helge Holden.
\newblock The {K}olmogorov-{R}iesz compactness theorem.
\newblock {\em Expo. Math.}, 28(4):385--394, 2010.

\bibitem{MR4383118}
Tuomas Hyt\"{o}nen and Stefanos Lappas.
\newblock Extrapolation of compactness on weighted spaces: bilinear operators.
\newblock {\em Indag. Math. (N.S.)}, 33(2):397--420, 2022.

\bibitem{HytonenOikariSinko}
Tuomas Hyt\"onen, Tuomas Oikari, and Jaakko Sinko.
\newblock Fractional {B}loom boundedness and compactness of commutators.
\newblock 2022.
\newblock arXiv:2207.01385.

\bibitem{IsPoTr}
Joshua Isralowitz, Sandra Pott, and Sergei Treil.
\newblock {Commutators in the two scalar and matrix weighted setting}.
\newblock {\em preprint}, 2017.
\newblock \textsf{arXiv:2001.11182}.

\bibitem{MR1187088}
Tadeusz Iwaniec.
\newblock {$L^p$}-theory of quasiregular mappings.
\newblock In {\em Quasiconformal space mappings}, volume 1508 of {\em Lecture
  Notes in Math.}, pages 39--64. Springer, Berlin, 1992.

\bibitem{MR1835563}
Steven~G. Krantz and Song-Ying Li.
\newblock Boundedness and compactness of integral operators on spaces of
  homogeneous type and applications. {I}.
\newblock {\em J. Math. Anal. Appl.}, 258(2):629--641, 2001.

\bibitem{Lern2012}
Andrei~K. Lerner.
\newblock On an estimate of {C}alder{\'o}n-{Z}ygmund operators by dyadic
  positive operators.
\newblock {\em J. Anal. Math.}, 121(1):141--161, 2013.

\bibitem{LORR}
Andrei~K. Lerner, Sheldy Ombrosi, and Israel~P. Rivera-R\'{\i}os.
\newblock On two weight estimates for iterated commutators.
\newblock {\em J. Funct. Anal.}, 281(8):Paper No. 109153, 46, 2021.

\bibitem{perez94}
Carlos P{\'e}rez.
\newblock Two weighted inequalities for potential and fractional type maximal
  operators.
\newblock {\em Indiana Univ. Math. J.}, 43(2):663--683, 1994.

\bibitem{MR1481632}
Carlos P\'{e}rez.
\newblock Sharp estimates for commutators of singular integrals via iterations
  of the {H}ardy-{L}ittlewood maximal function.
\newblock {\em J. Fourier Anal. Appl.}, 3(6):743--756, 1997.

\bibitem{MR467384}
Akihito Uchiyama.
\newblock On the compactness of operators of {H}ankel type.
\newblock {\em Tohoku Math. J. (2)}, 30(1):163--171, 1978.

\bibitem{MR3834654}
Huoxiong Wu and Dongyong Yang.
\newblock Characterizations of weighted compactness of commutators via {${\rm
  CMO}(\Bbb R^n)$}.
\newblock {\em Proc. Amer. Math. Soc.}, 146(10):4239--4254, 2018.

\bibitem{MR4303944}
Qingying Xue, K\^{o}z\^{o} Yabuta, and Jingquan Yan.
\newblock Weighted {F}r\'{e}chet-{K}olmogorov theorem and compactness of
  vector-valued multilinear operators.
\newblock {\em J. Geom. Anal.}, 31(10):9891--9914, 2021.

\end{thebibliography}

\end{document}